\documentclass[12pt]{article} 
\usepackage{fourier}
\usepackage{index}
\usepackage[expansion=false]{microtype}
\usepackage[reqno]{amsmath}
\usepackage{fixmath}
\usepackage{amssymb, amsthm, enumerate}
\usepackage{mathtools}

\usepackage{hyperref}

\newtheoremstyle{plainsl}%
	{\topsep}
	{\topsep}
	{\slshape} 
	{}
	{\normalfont\bfseries}
	{.}
	{ }
	{}

\swapnumbers

{\theoremstyle{plainsl}
\newtheorem{theorem}{Theorem}[section]
\newtheorem{lemma}[theorem]{Lemma}
\newtheorem{corollary}[theorem]{Corollary}}
{\theoremstyle{remark}
}

\newcommand\lref[1]{Lemma~\ref{lem:#1}}
\newcommand\tref[1]{Theorem~\ref{thm:#1}}
\newcommand\cref[1]{Corollary~\ref{cor:#1}}

\renewcommand\proof{\noindent\textsl{Proof. }}
\newcommand\sqr[2]{{\vbox{\hrule height.#2pt
    \hbox{\vrule width.#2pt height#1pt \kern#1pt
        \vrule width.#2pt}\hrule height.#2pt}}}
\renewcommand\qed{%
	\ifmmode\eqno\sqr53
	\else\nolinebreak\ \hfill\sqr53\medbreak\fi}

\numberwithin{equation}{section}

\newcommand\al{\alpha}
\newcommand\be{\beta}
\newcommand\de{\delta}
\newcommand\De{\Delta}

\newcommand\ga{\gamma}

\newcommand\sg{\sigma}

\renewcommand\th{\theta} 

\newcommand\vphi{\varphi}

\newcommand\cx{{\mathbb C}}

\newcommand\ints{{\mathbb Z}}
\newcommand\re{{\mathbb R}}
\newcommand\rats{{\mathbb Q}}
\newcommand\comp[1]{{\mkern2mu\overline{\mkern-2mu#1}}}
\newcommand\diff{\mathbin{\mkern-1.5mu\setminus\mkern-1.5mu}}
\newcommand\sbs{\subseteq}

\newcommand\seq[3]{#1_{#2},\ldots,#1_{#3}}

\DeclareMathOperator{\aut}{Aut}

\newcommand\pmat[1]{\begin{pmatrix} #1 \end{pmatrix}}
\newcommand\sm[3]{\sum_{#1=#2}^{#3}}
\DeclareMathOperator{\rk}{rk}
\DeclareMathOperator{\tr}{tr}

\DeclareMathOperator{\sym}{Sym}
\DeclareMathOperator{\dist}{dist}

\newcommand\one{{\bf1}}

\newcommand\cprod{\mathbin{\square}}
\newcommand\vc[1]{|#1\rangle}

\title{State Transfer on Graphs} 

\author{
	Chris Godsil\\
	Combinatorics \& Optimization\\
	University of Waterloo}

\begin{document}
\maketitle
	
\begin{abstract}
    If $X$ is a graph with adjacency matrix $A$, then we define $H(t)$ to be
    the operator $\exp(itA)$. We say that we have perfect state transfer in
    $X$ from the vertex $u$ to the vertex $v$ at time $\tau$ if the $uv$-entry
    of $|H(\tau)_{u,v}|=1$. State transfer has been applied to key distribution
    in commercial cryptosystems, and it seems likely that other applications will be found.
	We offer a survey of some of the work on perfect state transfer and related questions.
	The emphasis is almost entirely on the mathematics.
\end{abstract}

\section{Perfect State Transfer}

Let $X$ be a graph on $n$ vertices with with adjacency matrix $A$ and 
let $H(t)$ denote the matrix-valued function $\exp(iAt)$. If $u$ and $v$ are distinct
vertices in $X$, we say \textsl{perfect state transfer} from $u$ to $v$
occurs if there is a time $\tau$ such that $|H(\tau)_{u,v}|=1$. (We will occasionally
use ``pst'' as an abbreviation for ``perfect state transfer''.)
We say that $X$ is \textsl{periodic} relative to a vertex $u$
if there is a time $\tau$ such that $|H(\tau)_{u,u}|=1$, and we
say $X$ itself is periodic if there is a time $\tau$ such that
$|H(\tau)_{u,u}|=1$ for all vertices $u$. 

We can use the complete graph $K_2$ as an illustration. Here
\[
    A = \pmat{0&1\\ 1&0},
\]
so $A^n$ equals $I$ if $n$ is even and equals $A$ if $n$ is odd. Consequently
\[
    H(t) = \cos(t)I + i\sin(t)A = \pmat{\cos(t)& i\sin(t)\\ i\sin(t)& \cos(t)}
\]
and hence
\[
    H(\pi/2) = \pmat{0&i\\ i&0}.
\]
This shows that we have perfect state transfer from $u$ to $v$ at time $\pi/2$.
We also see that $X$ is periodic with period $\pi$ (because $H(\pi)=-I$).
Of course at $t=\pi/2$ we also have perfect state transfer from $v$ to $u$; we will
see that this is not an accident.

We note two properties of $H(t)$:
\begin{enumerate}[(a)]
    \item
    Since $A$ is symmetric, $H(t)$ is symmetric.
    \item
    Since $\comp{\exp(itA)}=\exp(-itA)$, we find that $H(t)$ is unitary.
\end{enumerate}

\begin{lemma}
    If we have perfect state transfer on $X$ from $u$ to $v$ at time $\tau$, then
    we have perfect state transfer from $v$ to $u$ at the same time, and $X$ is periodic
    at $u$ and $v$ with period dividing $2\tau$.
\end{lemma}

\proof
If $u\in V(X)$, let $\vc{u}$ denote the vector that is one on $u$ and zero elsewhere.
(So $\vc{u}$ is the characteristic vector of $u$, viewed as a subset of $V(X)$.)
If we have pst from $u$ to $v$ at time $\tau$ then there is a complex number $\ga$ 
of norm 1 such that
\[
    H(\tau)\vc{u} = \ga\vc{v}.
\]
But this says that $H(\tau)_{u,v}=\ga$ and, since $H(\tau)$ is symmetric, we have
$H(\tau)_{v,u}=\ga$ and therefore there is pst from $v$ to $u$ at time $\tau$.

Now we see that
\[
    H(\tau)^2\vc{u} = \ga^2\vc{u}, \qquad H(\tau)^2\vc{v} = \ga^2\vc{v}
\]
and since $H(\tau)^2=H(2\tau)$, we've proved our second claim.\qed

If $|H(t)_{u,v}|=1$ then since $H(t)$ is unitary, the $uv$-entry of $H(t)$ is the only
non-zero entry in the $u$-column and the only non-zero entry in the $v$ row. Hence
if we have perfect state transfer from $1$ to $2$ at time $\tau$, then $H(\tau)$
has the form
\[
    H(\tau) = \pmat{T&0\\ 0&H_1}
\]
where $H_1$ is unitary and
\[
    T = \ga \pmat{0&1\\ 1&0}
\]
with $|\ga|=1$.

The theory of perfect state transfer starts with the papers Bose \cite{Bose2003}
and Christandl et al~\cite{Christandl2005}; we will refer to the latter frequently. 

We note that our matrix $H(t)$ determines what is known as a \textsl{continuous
quantum walk}, for background on this we refer the reader to \cite{Kempe2003, Kendon2003}.
Physicists use $\exp(-itA)$ where we have used $\exp(itA)$; this makes absolutely
no difference to the theory, which is what we care about here. 
For recent surveys on state transfer see Kendon and Tamon \cite{KendTam2011}, 
Stevanovi\'c \cite{Stevanovic2011} and Kay \cite{Kay2010}.

\section{Spectral Decomposition}

The main tool we use is spectral decomposition of symmetric matrices. Suppose $A$
is symmetric with distinct eigenvalues
\[
    \seq\th1m
\]
and let $E_r$ denote orthogonal projection on the eigenspace belonging to $\th_r$.
Then
\[
    E_r^2 = E_r = E_r^\dagger
\]
and if $f$ is a complex-valued function defined on the eigenvalues of $A$,
\[
    f(A) = \sm r1m f(\th_r)E_r.
\]
Taking $f$ to be the exponential matrix we obtain the basic identity
\[
    H(t) = \sm r1m \exp(i\th_r t)E_r.
\]
One important consequence of this is that, for each $t$, the matrix $H(t)$ is a
polynomial in $A$. Thus it commutes with $A$ and, more generally, with any matrix
that commutes with $A$. Further
\[
    H(t)\vc{u} = \sm r1m \exp(i\th_r t)E_r\vc{u}
\]
where the non-zero vectors $E_r\vc{u}$ are eigenvectors for $A$ (and $H$). The
set of eigenvalues $\th_r$ such that $E_r\vc{u}\ne0$ is the
\textsl{eigenvalue support} of the vector $\vc{u}$.

\begin{lemma}\label{lem:prjs}
    Let $\seq E1m$ be the idempotents in the spectral decomposition of $A(X)$
    and let $\seq\th1m$ be the corresponding eigenvalues.
    Then there is perfect state transfer from $u$ to $v$ at time $\tau$ if and
    only if there is a constant $\ga$ such that 
    \[
        E_r\vc{u} = \ga\exp(-i\tau\th_r) E_r\vc{v} \qquad (r=1,\ldots,m)
    \]
\end{lemma}

\proof
We have
\[
    H(\tau)\vc{u} = \ga\vc{v}.
\]
if and only if, for all $r$,
\[
    \ga E_r\vc{v} = E_rH(\tau)\vc{u}.
\]
Since 
\[
    E_rH(\tau) = H(\tau)E_r = \exp(i\th_r\tau)E_r
\]
the lemma follows.\qed

\begin{corollary}\label{cor:eqsup}
    If there is perfect state transfer from $u$ to $v$, then 
    $E_r\vc{u}=\pm E_r\vc{v}$, and accordingly $u$ and $v$ have the same
    eigenvalue support.
\end{corollary}

\proof
If we have state transfer from $u$ to $v$, then $E_r\vc{u}=\be E_r\vc{v}$
where $|\be|=1$. As $E_r\vc{u}$ and $E_r\vc{v}$ are both real, $\be=\pm1$.\qed

For later use we note some properties of the eigenvalue support of a vertex, but
to prove these we will need to provide an expression for the idempotents $E_r$.
If $\seq\th1m$ are the distinct eigenvalues of $A$, define the polynomial $p_k(t)$
by
\[
    p_k(t) = \prod_{r\ne k} \frac{t-\th_r}{\th_k-\th_r}.
\]
Then it is not hard to verify that $E_r=p_r(A)$.

\begin{lemma}\label{lem:esupp}
    Suppose $u\in V(X)$ and $S$ is its eigenvalue support. If $\th\in S$ then
    all algebraic conjugates of $\th$ are in $S$. If $X$ is bipartite and $\th\in S$
    then $-\th\in S$. The spectral radius of the connected component of $X$ that contains
    $u$ belongs to $S$.
\end{lemma}

\proof
If $X$ is not connected, then the elements of $S$ are eigenvalues of the connected 
component of $X$ that contains $u$, and the associated eigenvectors are zero on 
vertices not in this component. So we may assume $X$ is connected.
If $\th_r$ and $\th_s$ are algebraic conjugates then $E_r=p_r(A)$ and $E_s=p_s(A)$ 
are algebraic
conjugates and so $E_r\vc{u}\ne0$ if and only if $E_s\vc{u}\ne0$. The eigenvalue
belonging to the spectral radius of a connected graph is simple and the corresponding
eigenvector has no entry zero. It follows that no entry of the associated idempotent
is zero, which implies the third claim.

Suppose $X$ is bipartite on $n$ vertices. Let $D$ be the $n\times n$ diagonal
matrix such that $D_{v,v}=1$ if $v$ is at even distance from $u$, and $D_{v,v}=-1$
otherwise. If $DAD=-A$ and if $Az=\th z$ then $ADz = -\th Dz$.\qed

\section{Period}

If we have perfect state transfer from $u$ to $v$ in $X$ at time $\tau$, then
$X$ is periodic at $u$ with period $2\tau$. Our next result shows that minimum time
at which perfect state transfer involving $u$ occurs is determined by the minimum
period at $u$.

Some preliminaries. Assume $X$ is connected and
let $T$ denote the set of times $\tau$ such that $H(\tau)e_u$ is a scalar multiple
of $e_u$. Then $T$ is an additive subgroup of $\re$. Since for small $t$,
\[
    H(t) \approx I + itA
\]
and since $Ae_u\ne0$, it follows that $T$ is a discrete subgroup of $\re$. Hence it
is cyclic, generated by the minimum period of $X$ at $u$. If we have perfect state
transfer from $u$ to $v$ at time $\tau$ then $X$ is periodic at $u$ with period $2\tau$,
and hence $\tau\ge\sg/2$.

\begin{lemma}
    Suppose $X$ is a connected graph and $X$ is periodic at $u$ with minimum period $\sg$.
    Then if there is perfect state transfer from $u$ to $v$, there is perfect state
    transfer from $u$ to $v$ at time $\sg/2$.
\end{lemma}

\proof
Suppose we have $uv$-pst with minimum time $\tau$. Then $X$ is periodic at $u$,
with minimum period $\sg$ (say).

If $\sg<\tau$, then $H(\tau-\sg)e_u=\ga e_v$ for some $\ga$ and so $\tau$ is not
minimal. Hence $\tau<\sg$. Since $X$ is periodic at $\tau$ with period $2\tau$,
we see that $\sg\le 2\tau$. If $\sg<2\tau$ then $u$ is periodic with period dividing 
$2\tau-\sg$ and so $\sg\le2\tau-\sg$, which implies that $\sg\le\tau$. We conclude 
that $\sg=2\tau$.

Thus if the minimum period of $X$ at $u$ is $\sg$ and there is perfect state transfer
from $u$ to $v$, then there is perfect state transfer from $u$ to $v$ at time $\sg/2$
(and not at any shorter time).\qed

We have the following corollary, due to Kay \cite[Section~IIID]{Kay2011}:

\begin{corollary}
    If we have perfect state transfer in $X$ from $u$ to $v$ and also from $u$ to $w$,
    then $v=w$.\qed
\end{corollary}

It is actually possible to derive a lower bound on the minimum period in terms of the
eigenvalues of $X$.

\begin{lemma}
    \label{lem:minper}
    If $X$ is a graph with eigenvalues $\seq\th1m$ and transition matrix $H(t)$.
    If $x$ is a non-zero vector, then the minimum time $\tau$ such that 
    $x^TH(\tau)x=0$ is at least $\frac\pi{\th_1-\th_m}$.
\end{lemma}

\proof
Assume $\|x\|=1$. We want
\[
    0 = x^T H(t)x = \sum e^{it\th_s} x^TE_s x,
\]
where the sum is over the eigenvalues $\th_s$ such that $E_rx\ne=0$, i.e., over
the eigenvalue support of $x$. Since
\[
    1 = x^Tx = \sum x^TE_s x,
\] 
the right side is a convex combination of complex numbers of norm 1. 
When $t=0$ these numbers are all equal to 1, and as $t$ increases they spread out on 
the unit circle of radius. If they are contained in an arc of length
less than $\pi$, their convex hull cannot contain $0$, and for small(ish) values
of $t$, they lie in the interval bounded by $t\th_1$ and $t\th_m$. So for $x^TH_X(t)x$ 
to be zero, we need $t(\th_1-\th_m)\ge\pi$, and thus we have the constraint
\[
    t \ge \frac\pi{\th_1-\th_m}.\qed
\]

If $u\in V(X)$ and $x=\vc{u}$, then this bound is tight for $P_2$ but not for $P_3$.

\begin{lemma}
    If $X$ is a graph with eigenvalues $\seq\th1m$, the minimum period of $X$ at a vertex
    is at least $\frac{2\pi}{\th_1-\th_m}$.
\end{lemma}

\proof
We want
\[
    \ga = \sum_r e^{i\th_r t}(E_r)_{u,u},
\]
where $\|\ga\|=1$, and for this to hold there must be integers $m_{r,s}$ such that
\[
    t(\th_r-\th_s) = 2m_{r,s}\pi.
\]
This yields the stated bound.\qed

In the previous lemma, $\th_1$ is the spectral radius of $A(X)$. However $\th_m$
can be replaced by the least eigenvalue in the eigenvalue support. If the entries of $x$
are non-negative, these comments apply to \lref{minper} too. For more bounds along
the lines of the last two lemmas, go to \cite[Section~IIIC]{Kay2011}.

\section{More Examples}

Our theory is developing nicely, but as yet we have just one example of perfect state
transfer. We describe a second, also from Christandl et al.~\cite{Christandl2005}.

\begin{lemma}
    There is perfect state transfer between the end vertices of the path on three
    vertices at time $\pi/\sqrt2$.
\end{lemma}

\proof
The eigenvalues of $P_3$ are $\sqrt2$, $0$, $-\sqrt2$ with respective eigenvectors
\[
    \frac12\pmat{1\\ \sqrt2\\1},\quad \frac1{\sqrt2}\pmat{1\\0\\-1},
        \quad \frac12\pmat{1\\-\sqrt2\\1}.
\]
If we denote these vectors by $z_1$, $z_2$, $z_3$ respectively, then
\[
    E_r = z_rz_r^T.
\]
Then
\[
    H(t) =\exp(it\sqrt2)E_1 +E_2 +\exp(-it\sqrt2)E_3
\]
and consequently
\[
    H(\pi/\sqrt2) = -E_1+E_2-E_3 = \pmat{0&0&-1\\ 0&-1&0\\ -1&0&0}.\qed
\]

If $X$ and $Y$ are graphs then their \textsl{Cartesian product} $X\cprod Y$ is defined
as follows. Its vertex set is $V(X)\times V(Y)$, and $(x_1,y_1)$ is adjacent to $(x_2,y_2)$
if either
\begin{enumerate}[(a)]
    \item 
    $x_1=x_2$ and $y_1$ is adjacent to $y_2$, or
    \item
    $x_1$ is adjacent to $x_2$ and $y_1=y_2$.
\end{enumerate}

Thus the Cartesian product of the paths $P_m$ and $P_n$ is the $m\times n$ grid.
We use $X^{\square d}$ to denote the $d$-th Cartesian power of $X$---the Cartesian product
of $d$ copies of $X$.
The $d$-th Cartesian power of $P_2$ is the $d$-cube $Q_d$.

The theory of the Cartesian product is very well developed, 
for details for \cite{Imrich2000}.
We could have defined the Cartesian product using adjacency matrices:
\[
    A(X\cprod Y) = A(X)\otimes I + I\otimes A(Y).
\]
This expresses $A(X\cprod Y)$ as the sum of two commuting matrices, and hence allows
one to prove the following (due once again to Christandl et al.).

\begin{lemma}
    For any graphs $X$ and $Y$ we have
    \[
        H_{A(X\cprod Y)}(t) = H_{A(X)}(t) \otimes H_{A(Y)}(t).\qed
    \]
\end{lemma}

This lemma is particularly important in physical terms, because it implies that
a physical system modelled by $X\cprod Y$ is a composite of the systems modelled
by $X$ and $Y$. It also means that we now have infinitely many examples where 
perfect state transfer occurs. The vertices of the $d$-th Cartesian power
of $X$ are the $d$-tuples in $V(X)^d$, and if $u$ and $v$ are two such $d$-tuples,
then the distance between them is
\[
    \sum_{r=1}^d \dist_X(u_r,v_r).
\]

\begin{theorem}
    If we have perfect state transfer from $u$ to $v$ in $X$ at time $\tau$,
    then at time $\tau$ we have perfect state transfer in the $d$-th Cartesian 
    power of $X$ between the $d$-tuples
    \[
        (u,\ldots,u),\quad (v,\ldots,v).\qed
    \]
\end{theorem}

Since we have perfect state transfer on $P_2$ and $P_3$, we might naturally expect
that perfect state transfer is possible on all paths. We will see that this
is false.

\section{Periodicity}

We have seen that the existence of perfect state transfer implies periodicity.
Studying periodicity is an effective stepping stone to the study of state transfer,
and so we take this step.

The first thing to note is that
\[
    \sum E_r = I
\]
and so $X$ itself will be periodic if there is a time $\tau$ and a scalar $\ga$
with $|\ga|=1$ such that
\[
    e^{i\tau\th_r} = \ga,\quad r=1,\ldots,m.
\]
Certainly taking $\tau$ to be zero works. What is much more interesting is that if the 
eigenvalues of $X$ are integers then $\tau=2\pi$ works: if the eigenvalues of $X$ are integers
then $X$ is periodic with period dividing $2\pi$. The path $P_2=K_2$ is an example.

With only a little more thought we see that if there there is a number $\de$
such that $\th_r/\de\in\rats$ for all $r$, then $X$ is periodic with period dividing
$2\pi/\de$. The basic question is to what extent this rationality condition is
necessary.

The ratio condition is a necessary condition for a graph to be periodic
at a vertex. The version we offer here is stated as Theorem~2.2 in \cite{Godsil2008};
it is an extension of result from Saxena, Severini and Shparlinski \cite{Saxena2007},
which in turn extends an idea used in Christandl et al.~\cite{Christandl2005}.

\begin{theorem}\label{thm:rat}
    Let $X$ be a graph and let $u$ be a vertex in $X$ at which $X$ is periodic.
    If $\th_k$, $\th_\ell$, $\th_r$, $\th_s$ are eigenvalues in the support
    of $\vc{u}$ and $\th_r\ne\th_s$, then
    \[
        \frac{\th_k-\th_\ell}{\th_r-\th_s} \in \rats.\qed
    \]
\end{theorem}

Using this one can prove that a graph is periodic if and only if the ratio of any
two eigenvalues is rational, and this leads to the following result from \cite{Godsil2008}

\begin{theorem}
    A graph $X$ is periodic if and only if either:
    \begin{enumerate}[(a)]
        \item 
        The eigenvalues of $X$ are integers, or
        \item
        The eigenvalues of $X$ are rational multiples
        of $\sqrt\De$, for some square-free integer $\De$.
    \end{enumerate}
    If the second alternative holds, $X$ is bipartite.\qed
\end{theorem}

If $X$ is regular then its spectral radius is an integer and so (b) cannot hold.
Thus a regular graph is periodic if and only if its eigenvalues are integers.

Since our actual concern is perfect state transfer, not periodicity, it
will only be useful if there are interesting cases where perfect state transfer
implies periodicity (and just periodicity at a vertex). We take this up in the next
section. Note that there are graphs with perfect state transfer that
are not periodic. Stevanovi\'c observes that the bipartite complements
of an even number of copies of $P_3$ provide a family of examples, and
another class is presented in Angeles-Canul et al.~\cite{Angeles-Canul2009}.
(A bipartite graph $Y$ is a \textsl{bipartite complement} of a bipartite graph $X$ 
with bipartition $(A,B)$ if the edge set $E(Y)$ is the complement
of $E(X)$ in the edge set of the complete bipartite graph with bipartition $(A,B)$.)

\section{Vertex-Transitive Graphs}

An automorphism of the graph $X$ is a permutation $\al$ of $V(X)$ such that
the vertices $u^\al$ and $v^\al$ are adjacent if and only if $u$ and $v$ are.
Any permutation can be represented by a permutation matrix, and a permutation matrix
$P$ is an automorphism of $X$ if and only if it commutes with $A(X)$. The set of
all automorphisms of $X$ forms its automorphism group $\aut(X)$. Our graph $X$ is 
\textsl{vertex transitive} if $\aut(X)$ is transitive as a permutation group, that is,
for each pair of vertices $u$ and $v$ there is an automorphism $\al$ such that $u^\al=v$.

\textsl{Cayley graphs} form an important class of vertex-transitive graphs.
To construct a Cayley graph for a group $G$ we first choose a subset $C$ of $G$.
The vertex set of the Cayley graph $X(G,C)$ is $G$, and elements $g$ and $h$ of $G$
are adjacent if $hg^{-1}\in C$. We call $C$ the \textsl{connection set}, and we do
\textbf{not} assume that $C$ generates $G$ (so $X$ might not be connected).
To avoid loops and multiple edges we do assume that $1\notin C$ and the $C$ is
inverse-closed; if $g\in C$ then $g^{-1}\in C$. If $a\in G$ then the map
that send $g$ to $ga$ is an automorphism of $G$. In fact $G$ acts regularly on $V(X)$
by right multiplication, and so any Cayley graph is vertex transitive.

There are two classes of Cayley graphs which are important to us. If $G$ is the cyclic
group $\ints_n$, then a Cayley graph for $G$ is a \textsl{circulant}. If $G$ is
the elementary abelian 2-group $\ints_2^d$, then $X$ is a so-called
\textsl{cubelike graph}. The cycle on $n$ vertices is a circulant with 
connection set $\{1,-1\}$,
and the $d$-cube is a cubelike graph. Note that if $G$ is abelian we are using $+$
as our group operation.

For vertex-transitive graph, the existence of perfect state transfer
has very strong consequences, as shown by the following result. (This
is a consequence of \cite[Theorem~4.1]{Godsil2008}.)

\begin{theorem}
    \label{thm:hcf2}
	Suppose $X$ is a connected vertex-transitive graph with vertices $u$ and $v$, 
	and perfect
	state transfer from $u$ to $v$ occurs at time $\tau$. Then $H(\tau)$ is a 
	scalar multiple of a permutation matrix with order two and no fixed points,
	and it lies in the centre of the automorphism group of $X$.\qed
\end{theorem}

An immediate consequence is that if perfect state transfer takes place on a 
vertex-transitive graph $X$, then $|V(X)|$ is even. 

We can weaken the assumption that $X$ is vertex transitive in this theorem
it is enough that $A(X)$ should belong to a homogeneous coherent algebra.
A \textsl{coherent algebra} is a vector space of matrices that is closed under
both the usual matrix multiplication and under Schur multiplication and contains
$I$ and $J$. Such an algebra has a unique basis of $01$-matrices and it is
\textsl{homogeneous} if $I$ is an element of this basis (rather than a sum of elements).
The adjacency matrix of a vertex-transitive graph belongs to a homogeneous
coherent algebra, and so does the adjacency matrix of a distance-regular graph. 
For details see \cite[Theorem~4.1]{Godsil2008}.

For vertex-transitive graphs we can specify the true period: if the eigenvalues 
are integers and $2^e$ is the largest power of $2$ that divides the greatest 
common divisor of the eigenvalues then the period is $\pi/2^{d-1}$.

\section{Cayley Graphs of Abelian Groups}

In investigations of state transfer, Cayley graphs of abelian groups provide
a useful test bed, because it is easy to compute their eigenvalues and eigenvectors.

Another advantage of this class is that we can decide which vertices might be involved
in state transfer. Each finite group $G$ gives rise to two regular permutation groups: 
the group of permutations given by right multiplication on $G$ and the group of 
permutations given
by left multiplication. If $n=|G|$, these give two regular subgroups of the symmetric
group $\sym(n)$ and each element in one group commutes with each element in the other. 
The intersection of these two groups consists of the elements in the center of $G$. 
(So if $G$ is abelian, the two groups are equal.) 
If $T$ denotes the permutation matrix arising in \tref{hcf2}, then we have
the following extension of this theorem:

\begin{lemma}
    If $X$ is a Cayley graph for a group $G$ and perfect state transfer occurs
    at time $\tau$, then $T$ lies in the center of $G$.\qed
\end{lemma}

Even when $G$ is abelian this is useful, because it tells us that $T$ is an element of $G$.
So if we have perfect state transfer on a Cayley graph for an abelian group $G$, then
it maps the vertex $g$ to $g+c$ for some element of order two in $Z(G)$. If $G$ is the
cyclic group of order $n=2\nu$, then perfect state transfer must send $0$ to $\nu$
(and $a$ to $a+\nu$). If $G$ is cyclic then this element of order two is unique
(it if it exists).

\subsection{Cubelike Graphs}

A cubelike graph is a Cayley graph for $\ints_2^d$. The adjacency matrix of
a cubelike graph can written as sum of commuting permutation matrices $P$ such 
that $P^2=I$ and $\tr(P)=0$. If 
\[
    A =P_1+\cdots+P_d
\]
then
\[
    H(t) = \exp(it(P_1+\cdots+P_d)) =\prod_{r=1}^d \exp(itP_r).
\]
But if $P^2=I$ then
\[
    \exp(itP) = \cos(t)I + i\sin(t)P
\]
and therefore
\[
    H(t) = \prod_{r=1}^d (\cos(t)I + i\sin(t)P_r).
\]
Consequently $H(\pi)=(-1)^dI$ and
\[
    H(\pi/2) = i^d \prod_{r=1}^d P_r.
\]
Using these ideas we arrive at the following result from \cite{Bernasconi2008}

\begin{lemma}
    Suppose $C\sbs\ints_2^d\diff0$ and $X=X(\ints_2^d,C)$. Then $X$ is periodic
    with period dividing $\pi$. Its period is equal to $\pi$ if and only the
    sum $\sg$ of the elements of $C$ is not zero, and in this case we have perfect
    state transfer from $0$ to $\sg$ at time $\pi/2$.\qed
\end{lemma}

The connection set of a cubelike graph on $2^d$ vertices with valency $m$ can be
presented as a $d\times m$ matrix over $\ints_2$ with distinct columns;
the columns of the matrix are the connection set. If $M$ is such a matrix then 
its row space is a binary code. A binary code is \textsl{even} if the Hamming weight
of any code word is even and it is easy to show that the code is even if and
only if $M\one=0$, that, is, the columns of $M$ sum to zero. Hence we have perfect 
state transfer on a cubelike graph at $\pi/2$ if and only if its code is not even.

We may have perfect state transfer when the code is even. A binary code is
\textsl{doubly even} if all code words have weight divisible by four.
The code is self-orthogonal if $MM^T=0$ and again it is not hard to show
that a self-orthogonal code is doubly even if and only if the weight of
each row of $M$ is doubly even. In \cite{Cheung2010} it is proved that
if the code of a cubelike graph is self-orthogonal and even, but not doubly
even, then we have perfect state transfer at time $\pi/4$.

\subsection{Circulants}

A circulant is a Cayley graph for the cyclic group $\ints_n$.
As we saw above there is no perfect state transfer if $n$ is odd
and, if $n=2\nu$ then any state transfer must be from $a$ to $a+\nu$.

However we have got ahead of ourselves---a circulant is periodic if and only its
eigenvalues are integers. So we need to determine when this holds. Computing the
eigenvalues of a circulant is easy because we know its eigenvectors. Choose
a primitive $n$-th root of unity in $\cx$, say $\zeta$, and an integer $s$.
Then the function that maps $x$ in $\ints_n$ to $\zeta^{sx}$ is an eigenvector.
If the connection set is $C$, then the eigenvalue is
\[
    \sum_{x\in C} \zeta^{sx}.
\]

What is more surprising is that it is easy to characterize the connection
sets $C$ such that $X(C)$ has only integer eigenvalues. Define two elements
of an abelian group $G$ to be equivalent if they generate the same subgroup of $G$.
Bridges and Mena \cite{Bridges1982} showed that the eigenvalues of $X(G,C)$ are integers
if and only if $C$ is a union of equivalence classes. For cyclic groups, two elements 
are equivalent if and only if they have the same order.

In \cite{MR2645074, MR2561744, MR2523011, Petkovic2011} Ba\v{s}i\'c,  
Petkovi\'c and (in one case) Stevanovi\'c have investigated perfect state transfer 
on circulants. Among their many results, they proved that if $n$ is squarefree or 
is congruent to 2 modulo 4, there is no perfect state transfer.
More recently Ba\v{s}i\'{c} \cite{Basic2011} has completely characterized
the circulants on which perfect state transfer occurs.

A \textsl{bicirculant} is a graph on $2n$ vertices admitting an automorphism
of order $n$ with two orbits of length $n$. (So the Petersen graph is one example.)
Any circulant of even order is a bicirculant but in general a bicirculant graph
need not be vertex transitive. If $X$ is a bicirculant relative to an automorphism
$g$, then the subgraphs induced by the orbits of $g$ are circulants.
Angeles-Canul et al.~\cite{Angeles-Canul2009} define a circulant join to be a 
bicirculant where the two induced circulants are isomorphic, and give a condition
for such a graph to admit perfect state transfer.

\section{Equitable Partitions}

A partition $\pi$ of the vertices of a graph $X$ is \textsl{equitable}
if, for each ordered pair of cells $C_i$ and $C_j$ from $\pi$, all vertices 
in $C_i$ have the same number of neighbors in $C_j$. (So the subgraph of
$X$ induced by a cell is a regular graph, and the subgraph formed by the
vertices of two cells and the edges which join them is bipartite and
semiregular, that is, all vertices in the same color class have the same valency.) 
We note that the orbits of any group of automorphisms of $X$ 
form an equitable partition. The \textsl{discrete partition}, with each cell 
a singleton, is always equitable; the \textsl{trivial partition} with exactly one cell
is equitable if an only if $X$ is regular. For the basics concerning equitable 
partitions see for example \cite[Section~9.3]{cggrbk}. In particular the 
join of two equitable partitions is equitable, and consequently given any 
partition of $X$, there is a unique coarsest refinement of it---the join of 
all equitable partitions which refine it. (Note: here ``join'' refers to join in 
the lattices of partitions.)

If $\pi$ is a partition of $V(X)$, its \textsl{characteristic matrix} is the $01$-matrix
whose columns are the characteristic vectors of the cells of $\pi$, viewed
as subsets of $V(X)$. If we scale the columns of the characteristic matrix
so they are unit vectors, we obtain the \textsl{normalized characteristic matrix}
of $\pi$. If $P$ and $Q$ are respectively the characteristic and normalized
characteristic matrix of $\pi$, then $P$ and $Q$ have the same column space 
and $Q^TQ=I$. The matrix $QQ^T$ is block diagonal, and its diagonal blocks are
all of the form
\[
    \frac1r J_r
\]
where $J_r$ is the all-ones matrix of order $r\times r$, and the size 
of the $i$-th block is the size of the $i$-th cell of $\pi$.
The vertex $u$ forms a singleton cell of $\pi$ if and only if $Qe_u=e_u$.

We use $|\pi|$ to denote the number of cells of $\pi$.

\begin{lemma}
    Suppose $\pi$ is a partition of the vertices of the graph $X$,
    and that $Q$ is its normalized characteristic matrix. Then
    the following are equivalent:
    \begin{enumerate}[(a)]
        \item 
        $\pi$ is equitable.
        \item
        The column space of $Q$ is $A$-invariant.
        \item
        There is a matrix $B$ of order $|\pi|\times|\pi|$ such that $AQ=QB$.
        \item
        $A$ and $QQ^T$ commute.\qed
    \end{enumerate}
\end{lemma}

If $u\in V(G)$, we use $\De_u$ to denote the partition of the vertices by distance 
from $u$. The following result is proved in \cite{Godsil1987a}.

\begin{theorem}
    A graph $X$ is distance-regular if it is regular and for each vertex $u$ in $X$,
    the distance partition $\De_u$ is equitable.
\end{theorem}

(If $X$ is not regular but all distance partitions are equitable, then it is a
distance-biregular graph. For details see \cite{Godsil1987a}.)

Ge at al \cite{Ge2010} provide an application of the theory of equitable partitions
unrelated to what we consider here.

\section{Stabilizers}

If $u\in V(X)$, then $\aut(X)_u$ denotes the group of automorphisms of $X$ 
that fix $u$.  Our next result says that if perfect state transfer from $u$ 
to $v$ occurs, then any automorphism of $X$ that fixes $u$ must fix $v$ 
(and vice versa).

\begin{lemma}\label{lem:fixuv}
    Suppose perfect state transfer from $u$ to $v$ occurs on $X$ at time $\tau$.
    If $MA=AM$ and $M\vc{u}=\vc{u}$ then $M\vc{v}=\vc{v}$.
\end{lemma}

\proof
Since $M$ must commute with $E_r$,
\[
    ME_r\vc{u} = E_rM\vc{u} =E_r\vc{u}.
\]
By \lref{prjs} it follows that $ME_r\vc{v}=E_r\vc{v}$, and 
therefore $M\vc{v}=\vc{v}$.\qed

\begin{corollary}
	If $X$ admits perfect state transfer from $u$ to $v$, then 
	$\aut(X)_u=\aut(X)_v$.\qed
\end{corollary}

By way of example, suppose $X$ is a Cayley graph for an abelian
group $G$. The map that sends $g$ in $G$ to $g^{-1}$ is an automorphism
of $G$ and also an automorphism of $X$. The fixed points of this
automorphism are the elements of $G$ with order one or two.
So if perfect state transfer from 1 to a second vertex $a$ occurs,
then $a$ has order two. Consequently perfect state transfer cannot
occur on a Cayley graph for an abelian group of odd order. (Another
proof of this is presented as Corollary~4.2 in \cite{Godsil2008}.)

If there is perfect state transfer from $u$ to $v$, then it follows from 
the previous corollary that the orbit partitions of $\aut(X)_u$ and $\aut(X)_v$
are equal. Since equitable partitions can be viewed as a generalization of orbit 
partitions, it is not entirely unreasonable to view the following result as an 
extension of this fact.

\begin{corollary}
    Let $u$ and $v$ be vertices in $X$ and let $\pi_u$ and $\pi_v$ denote the coarsest
	equitable partition of $X$ in which $\{u\}$ (respectively $\{v\}$) is a cell of size one.
	If there is perfect state transfer from $u$ to $v$, then $\De_u=\De_v$.
\end{corollary}

\proof
Let $Q$ be the normalized characteristic matrix of the partition
$\pi_u$, let $R=QQ^T$ and assume $H(t)=\exp(iAt)$. Then $H(t)$ is a polynomial in $A$
and so $RH(t)=H(t)R$. If we have perfect state transfer from $u$
to $v$ at time $\tau$, then
\[
    H(\tau)e_u = \ga e_v
\]
and accordingly
\[
    \ga Re_v = RH(\tau)e_u = H(\tau) Re_u = H(\tau)e_u = \ga e_v.
\]
So $Re_v=e_v$ and this implies that $\{v\}$ is a cell of $\pi_u$. By symmetry,
$\{u\}$ is a cell in $\pi_v$ and it follows that $\pi_u=\pi_v$\qed

By way of example, consider a distance-regular graph $X$ with diameter $d$.
If $u\in V(X)$, then $\De_u$ is the distance partition
with respect to $u$ with exactly $d+1$ cells, where the $i$-th
cell is the set of vertices at distance $i$ from $u$.
We conclude that if perfect state transfer occurs on $X$, then for each vertex
$u$ there is exactly one vertex vertex at distance $d$ from it.
In particular perfect state transfer does not occur on strongly regular graphs.
(This also follows from some observations in \cite[Section~4]{Godsil2008}.)

\begin{lemma}\label{lem:eq-pst}
    Suppose $C$ and $D$ are cells of an equitable partition $\pi$ of $X$ and that we 
    have perfect state transfer from $u$ in $C$ to $v$ in $D$ at time $\tau$. If
    $y$ is the characteristic vector of the cell of $\pi$ that contains $u$, 
    then $H(\tau)y$ is the scalar
    multiple of the characteristic vector of the cell that contains $v$.
    Further these two cells have the same size.
\end{lemma}

\proof
Let $Q$ be the normalized characteristic matrix of the partition. Then
\[
    H Q\vc{u} = Q H\vc{u} = \ga Q\vc{v}.
\]
If $w\in V(X)$ then $Q\vc{w}$ is the normalized characteristic vector of the cell
containing $w$ and so our first claim holds. Since $H$ is unitary, $z$ and $Hz$
have the same length.\qed

\section{Finiteness}

From \cite{Godsil2008} we know that if a graph is periodic, then the squares of 
its eigenvalues are integers and, if the graph is not bipartite, 
the eigenvalues themselves are integers. A variant of this
fact was derived in \cite{Godsil2010}, it implies that that if perfect state transfer occurs
on $X$, then the spectral radius of $X$ is an integer or a quadratic irrational.

\begin{theorem}
    \label{thm:intevs}
    Suppose $X$ is a connected graph with at least three vertices
    where perfect state transfer
    from $u$ to $v$ occurs at time $\tau$. Let $\de$ be the dimension of the
    $A$-invariant subspace generated by $e_u$. Then $\de\ge3$ and either
    all eigenvalues in the eigenvalue support of $u$ are integers, or
    they are all of the form $\frac12(a+b_r\sqrt\De$) 
    where $\De$ is a square-free integer and $a$ and $b_r$ are integers.\qed
\end{theorem}

\begin{corollary}
    There are only finitely many connected graphs with maximum valency at most $k$
    where perfect state transfer occurs.
\end{corollary}

\proof
Suppose $X$ is a connected graph where perfect state transfer
from $u$ to $v$ occurs at time $\tau$ and let $S$ be the eigenvalue support
of $u$. If the eigenvalues in $S$ are integers then $|S|\le2k+1$, and if they
are not integers then $|S|<(2k+1)\sqrt{2}$. So the dimension of the $A$-invariant
subspace of $\re^{V(X)}$ generated by $\vc{u}$ is at most $\lceil(2k+1)\sqrt{2}\rceil$, 
and this is also a bound on the maximum distance from $u$ of a vertex in $X$.
If $s\ge1$, the number of vertices at distance $s$ from $u$ is at most
$k(k-1)^{s-1}$, and the result follows.\qed

\begin{corollary}\label{cor:bipsq}
    Suppose $X$ is a bipartite graph with spectral radius $\th_1$. If perfect
    state transfer takes place on $X$, then $\th_1^2\in\ints$.
\end{corollary}

\proof
Suppose $X$ is bipartite are we have perfect state transfer from $u$ to $v$.
Let $S$ be the eigenvalue support of $u$. By \tref{intevs} there is a squarefree 
integer $\De$ and integers $a$ and $b_r$ such that each eigenvalue in $S$ has the form
\[
    \frac12(a+b_r\sqrt{\De}).
\]
We assume by way of contradiction that $a\ne0$. The spectral radius $\th_1$
is then $(a+b_r\sqrt{\De})/2$ and from \lref{esupp} we see that $\th_1$ and $-\th_1$
belong to $S$. This implies that $a=0$.\qed

\section{Cospectral Vertices}

We use $\phi(X,x)$ to denote the characteristic polynomial of $A(X)$.
Vertices $u$ and $v$ in the graph $X$ are \textsl{cospectral} if
\[
	\phi(X\diff u,t) = \phi(X\diff v,t).
\]
Of course two vertices that lie in the same orbit of $\aut(X)$ are cospectral,
but there are many examples of cospectral pairs of vertices where this does 
not hold. A graph is \textsl{walk regular} if any two of its vertices 
are cospectral. Any strongly regular graph is walk regular.

We note the following identities:
\[
    \frac{\phi(X\diff u,t)}{\phi(X,t)} = ((tI-A)^{-1})_{u,u}
        =\sum_r (t-\th_r)^{-1} (E_r)_{u,u}.
\]
(Here the first inequality is a consequence of Cramer's rule and the second is 
spectral decomposition.) Since
\[
    (E_r)_{u,u} = \| E_r\vc{u}\|^2
\]
we see that $u$ and $v$ are cospectral if and only the projections $E_r\vc{u}$
and $E_r\vc{v}$ have the same length for each $r$. So \lref{prjs} yields
immediately:

\begin{lemma}\label{lem:pst-cosp}
	If $X$ admits perfect state transfer from $u$ to $v$, then 
	$E_ru=\pm E_rv$ for all $r$, and $u$ and $v$ are cospectral.\qed
\end{lemma}

We note that if the eigenvalues of $A$ are simple and 
\[
    \| E_r\vc{u}\|^2 = \| E_r\vc{v}\|^2
\]
then necessarily $E_r\vc{u}=\pm E_r\vc{v}$.

Let $A$ be the adjacency matrix of the graph $X$ on $n$ vertices and suppose 
$S$ is a subset of $V(X)$ with characteristic vector $\vc{S}$. The \textsl{walk matrix}
of $S$ is the matrix with columns
\[
    \vc{S}, A\vc{S},\ldots, A^{n-1}\vc{S}.
\]
We say that the pair $(X,S)$ is \textsl{controllable} if this walk matrix is
invertible. If $u\in V(X)$ then the walk matrix of $u$ is just the walk matrix
relative to the subset $\{u\}$ of $V(X)$. From our discussion at the start of the
proof of \tref{intevs},
we have that the rank of the walk matrix relative to $u$ is equal to the
size of the eigenvalue support of $u$. Thus if $(X,u)$ is controllable,
the eigenvalues of $A$ must be distinct.

\begin{lemma}\label{lem:wutwu}
    Let $u$ and $v$ be vertices of $X$ with respective walk matrices $W_u$ 
    and $W_v$. Then $u$ and $v$ are cospectral if and only if $W_u^TW_u=W_v^TW_v$.\qed
\end{lemma}

One consequence of this is that if $u$ and $v$ are cospectral and
$(X,u)$ is controllable, then so is $(X,v)$. We also have the following
(from \cite{Godsil2010a}):

\begin{corollary}
    If $u$ is a vertex in $X$ with walk matrix $W_u$, then $\rk(W_u)$
    is equal to the number of poles of the rational function
    $\phi(X\diff u,x)/\phi(X,x)$. Hence $(X,u)$ is controllable if and only
    if $\phi(X\diff u,x)$ and $\phi(X,x)$ are coprime.
\end{corollary}

\proof
By the spectral decomposition,
\[
    ((xI-A)^{-1})_{u,u} = \sum_r \frac{1}{x-\th_r}(E_r)_{u,u}.\qed
\]

The characteristic polynomial of the path $P_n$ satisfies the recurrence
\[
    \phi(P_{n+1},x) = x\phi(P_n,x) - \phi(P_{n-1},x)
\]
and using this it easy to show that either end-vertex of a path is controllable.
On the other hand, if $X$ has an eigenvalue of multiplicity greater than one
then no vertex in $X$ is controllable. The following result comes from \cite{Godsil2010}.

\begin{theorem}\label{thm:nocontrol}
    Let $X$ be a connected graph on at least four vertices.
    If we have perfect state transfer between distinct vertices $u$ and 
    $v$ in $X$, then neither $u$ nor $v$ is controllable.\qed
\end{theorem}

For information on controllability, go to \cite{Godsil2010b}.

\section{Transfer without Exponentials}

If $u$ and $v$ are controllable and cospectral we cannot get perfect state transfer
between them. Our next result shows that if these conditions hold, there is nonetheless
a symmetric orthogonal matrix $Q$ which commutes with $A$ and maps $\vc{u}$ to $\vc{v}$.

\begin{lemma}
    Let $u$ and $v$ be vertices in $X$ with respective walk matrices $W_u$ and $W_v$.
    If $u$ and $v$ are controllable and $Q:=W_vW_u^{-1}$ then $Q$ is a polynomial in $A$.
    Further $Q$ is orthogonal if and only if $u$ and $v$ are cospectral.
\end{lemma}

\proof
Let $C_\phi$ denote the companion matrix of the characteristic polynomial of $A$. Then
\[
    AW_u =W_u C_\phi
\]
for any vertex $u$ in $X$. Hence if $u$ and $v$ are controllable,
\[
    W_u^{-1}AW_u = W_v^{-1}AW_v
\]
and from this we get that
\[
    AW_vW_u^{-1} = W_vW_u^{-1}A.
\]
Since $X$ has a controllable vertex its eigenvalues are all simple, and so any matrix
that commutes with $A$ is a polynomial in $A$. This proves the first claim. 

From \lref{wutwu}, the vertices $u$ and $v$ are cospectral
if and only if $W_u^TW_u=W_v^TW_v$, which is equivalent to
\[
    Q^{-T} = W_v^{-T}W_u^T = W_vW_u^{-1} = Q.\qed
\]

This lemma places us in an interesting position. If $u$ and $v$ are cospectral and
controllable, there is an orthogonal matrix $Q$ that commutes with $A$ such that
$Q\vc{u}=\vc{v}$ but, by the result of the previous section, $Q$ cannot be equal to
a scalar multiple of $H(t)$ for any $t$.

\section{Pretty Good State Transfer}

We say we have \textsl{pretty good state transfer} from $u$ to $v$ if
there is a sequence $\{t_k\}$ of real numbers and a scalar $\ga$ such that
\[
    \lim_{k\to\infty} H(t_k)\vc{u} = \ga\vc{v}.
\]

As an example consider $P_4$ with eigenvalues
\[
    \th_1=\frac12(\sqrt{5}+1),\quad \th_2=\frac12(\sqrt{5}-1),
		\quad \th_3=\frac12(-\sqrt{5}+1),\quad \th_4=\frac12(-\sqrt{5}-1).
\]
Then by straightforward computation
\[
    E_1-E_2+E_3-E_4 = \pmat{0&0&0&1\\ 0&0&1&0\\ 0&1&0&0\\ 1&0&0&0};
\]
denote the right side by $T$. (We could prove that $T=\sum_r(-1)^{r-1}E_r$ by verifying 
that if $z_r$ is an eigenvector with eigenvalue $\th_r$, then $Tz_r=(-1)^{r-1}z_r$.)
Now choose integers $a$ and $b$ so that
\[
    \frac{a}{b} \approx \frac{1+\sqrt{5}}{2}.
\]
Then $b\th_1\approx a$ and
\[
    b\th_2 =b(\th_1-1) \approx a-b,\quad b\th_3 \approx b-a,\quad b\th_4\approx -a.
\]
For example take $a=987$ and $b=610$ and set $\tau = 610\pi/2$. Then
the values of $\tau\th_r$ are (approximately)
\[
    987\pi/2,\quad 377\pi/2,\quad -377\pi/2,\quad -987\pi/2
\]
and these are congruent modulo $2\pi$ to
\[
    3\pi/2,\quad \pi/2,\quad 3\pi/2,\quad \pi/2.
\]
Hence 
\[
    H(305\pi) \approx -iT.
\]
(The approximation is accurate to five decimal places.)

In general we have to choose $a$ and $b$ so that $a\cong3$ and $b\cong2$ modulo 4.
If $f_n$ is the $n$-th Fibonacci number with $f_0=f_1=1$ then modulo 4
\[
    f_{4m+2} \cong 2,\quad f_{4m+3}\cong 3
\]
and the ratios $f_{n+1}/f_n$ are the standard continued fraction approximation 
to $(\sqrt{5}+1)/2$. We conclude that we have pretty good state transfer between 
the end-vertices of $P_4$. We leave the reader the exercise of verifying that
there is also pretty good state transfer between the end-vertices of $P_5$
(for which the eigenvalues are $0$, $\pm1$, $\pm\sqrt3$).
In the next section we see that we do \textbf{not} 
have perfect state transfer on $P_n$ when $n\ge4$. Also, to get a good approximation
to perfect state transfer on $P_5$ the numerical evidence is that $t$ must be very
large. This suggests that pretty good state transfer will not be a satisfactory
substitute for perfect state transfer in practice.

Dave Morris has noted the following in a private communication.

\begin{lemma}
    If we have pretty good state transfer from $u$ to $v$, then 
    $E_r\vc{u}=\pm E_r\vc{v}$ for each $r$.
\end{lemma}

\proof
By assumption there is a sequence $\{t_k\}$ of real numbers such that 
\[
    \lim_{k\to\infty} H(t_k)\vc{u} =\ga \vc{v}.
\]
Since the unit circle is compact there is a subsequence $\{t_\ell\}$ and a
complex number $\zeta$ such that $\exp(it_\ell)\to\zeta$. Now
\[
    \zeta E_r\vc{u} = \lim_{\ell\to\infty}\exp(it_\ell\th_r)E_r\vc{u}
        = \lim_{\ell\to\infty}E_r H(t_\ell)\vc{u} = \ga E_r\vc{v}.
\]
Since $E_r\vc{u}$ and $E_r\vc{v}$ are real vectors, the lemma follows.\qed

\section{Paths}

In \cite{Christandl2005} Christandl et al.~proved that perfect state transfer
between the end-vertices of a path on $n$ vertices did not occur if $n\ge4$.
It is possible to extend their arguments to show that $P_2$ and $P_3$
are the only paths where perfect state transfer occurs at all. Our approach
has benefited from discussions with Dragan Stevanovi\'c.

To begin we note some simple properties of paths. First, the characteristic
polynomials $\phi(P_n,x)$ satisfy the recurrence
\[
    \phi(P_{n+1},x) = x\phi(P_n,x) - \phi(P_{n-1},x).
\]
One consequence of this is $\phi(P_{n+1},x)$ and $\phi(P_n,x)$ are coprime.
Since interlacing implies that any multiple eigenvalue of $P_{n+1}$ must
be an eigenvalue of $P_n$, we also see that the eigenvalues of a path are simple.

Given the $\Phi(P_n,x)$ and $\phi(P_{n-1},x)$ are coprime, it follows
from \tref{nocontrol} that we do not have perfect state transfer between
end-vertices in $P_n$ when $n\ge4$. Using the identity
\[
    \phi(P_{m+n},x) = \phi(P_m,x)\phi(P_n,x)
        -\phi(P_{m-1},x)\phi(P_{n-1},x)
\]
(which can be derived easily by induction from our recurrence above), it is not
hard to show that $\phi(P_n,x)$ and $\phi(P_m,x)$ have a non-trivial common factor
if and only if $m+1$ divides $n+1$. This rules out many more possible cases of
perfect state transfer.

If $n$ is odd then the stabilizer in $\aut(P_n)$ of the middle vertex of $P_n$ 
has order two, while the stabilizer of any other vertex is trivial. 
So the middle vertex cannot be involved in perfect state transfer.

But nothing we have mentioned will rule out perfect state transfer between
(say) vertices 3 and 6 in $P_8$. We address this problem now.

Let $\zeta_n(x)$ denote the vector
\[
    \pmat{1\\ \phi(P_1,x)\\ \vdots\\ \phi(P_{n-1},x)}.
\]
If the eigenvalues of $P_n$ are
\[
    \th_1\ge \cdots \ge \th_n
\]
then $\zeta_n(\th_r)$ is an eigenvector for $P_n$ with eigenvalue $\th_r$.
This is easy to verify using the recurrence, and using this we can also see
that two consecutive entries of $\zeta_r$ cannot be zero. Nor can the last entry
be zero. We say there is a \textsl{sign change at $r$} in the sequence
$(\phi(P_r,x))_{r=0}^{n-1}$ if
\[
    \phi(P_{r-1}(x))\,\phi(P_r(x)) < 0
\]
or if $\phi(P_r(x))=0$ and
\[
    \phi(P_{r-1}(x))\,\phi(P_{r+1}(x)) < 0.
\]
From the recurrence we see that if $\phi(P_r(x))=0$ then there is a sign change at $r$.
It follows from Sturm's theorem that there are exactly $m-1$ sign changes in the sequence 
$(\phi(P_n,\th_m))_{r=0}^{n-1}$.

Let $T$ be the permutation matrix representing the non-identity automorphism of $P_n$.
Since the eigenvalues of the path are simple, if $z$ is an eigenvector for $A=A(P_n)$,
then $Tz$ is also an eigenvector and consequently $Tz = \pm z$. Since the
bottom entry of $\zeta_n(\th_2)$ must be negative, we have 
\[
    T\zeta_n(\th_2) = -\zeta_n(\th_2).
\]
We conclude that either no entry of $\zeta_n(\th_2)$ is zero, or $n$ is odd
and the middle entry is zero. Hence if $m$ is not the middle vertex
then neither $E_1\vc{m}$ nor $E_2\vc{m}$ is zero. If $D$ is the $n\times n$ diagonal
matrix with $D_{r,r}=(-1)^r$ and $z$ is an eigenvector of the path with eigenvalue $\th$,
then $Dz$ is an eigenvector with eigenvalue $-\th$. Hence $E_{n-1}\vc{m}$ and $E_n\vc{m}$
are not zero. By \tref{rat} we have that
\[
    \frac{\th_2}{\th_1} = \frac{\th_n-\th_{n-1}}{\th_1-\th_n} \in \rats.
\]
The eigenvalues of $P_n$ are the numbers
\[
    2\cos\left(\frac{r\pi}{n+1}\right),\qquad r=1,\ldots,n
\]
and it follows that
\[
    \th_2 = \th_1^2-2.
\]
Our rationality condition now implies that $\th_1$ must be, at worst, a quadratic
irrational. But by \cref{bipsq} we have that $\th_1^2$ is an integer, and therefore
$\th_2$ must also be an integer. Since $\th_2<2$ and $\th_1>0$ we have
\[
    \th_2 \in \{-1,0,1\}.
\]
If $\th_2=-1$ then $\th_1=1$ and $X=K_2$. If $\th_2=0$, then $\th_1=\sqrt2$
and $X=P_3$. If $\th_2=1$ then $\th_1=\sqrt{3}$ and according to the ratio
condition
\[
    \frac{\th_1-\th_2}{\th_1+\th_2} = \frac{\sqrt3-1}{\sqrt3+1} = 2-\sqrt3
\]
should be rational.

Depending on one's mood, it is either instructive or depressing to see how much effort 
is needed to deal with perfect state transfer on paths.

\section{Joins}

If $X$ and $Y$ are graphs let $X+Y$ denote their \textsl{join}, which
we get by taking a copy of $X$ and a copy of $Y$ and joining each vertex in $X$
to each vertex in $Y$. 
Angeles-Canul et al.~\cite{Angeles-Canul2009, Angeles-Canul2009a}
and the Ge et al.~\cite{Ge2010} provide many interesting results on perfect
state transfer in joins, including cases with weighted edges. 
Here we will focus simply on the joins of two regular graphs.

Assume $X$ is $k$-regular on $m$ vertices and $Y$ is $\ell$-regular on $n$ vertices.
Set $A=A(X)$ and $B=A(Y)$. If $Z:=X+Y$ then
\[
    A(Z) = \pmat{A&J\\ J^T&B}.
\]
If $Az=\th z$ and $\one^T z=0$, then
\[
    A(Z)\pmat{z\\0} = \th\pmat{z\\0}.
\]
Similarly, if $Bz=\th z$ then
\[
    A(Z)\pmat{0\\z} = \th\pmat{0\\z}.
\]
We see that $n+m-2$ of the eigenvalues of $X+Y$ are eigenvalues of $X$ and eigenvalues
of $Y$. The remaining two eigenvalues are associated with eigenvectors that are 
constant on $V(X)$ and $V(Y)$. Assume that $A$ and $B$ have respective spectral 
decompositions:
\[
    A = \sum_r \th_r E_r,\quad B = \sum_s \nu_s F_s
\]
where $E_1$ and $F_1$ are multiples of $J$. Then we have a decomposition
\[
    A(Z) =\mu_1 N_1 +\mu_2 N_2 +\sum_{r>1}\th_r\hat{E}_r +\sum_{s>1}\hat{F}_r,
\]
where a lot of explanation is needed. Here
\[
    \hat{E}_r = \pmat{E_r&0\\0&0},\quad \hat{F}_s = \pmat{0&0\\0&F_s}
\]
while since $\rk(N_1)=\rk(N_2)=1$ they are respectively of the form
\[
    \pmat{aJ_{m,m}&\sqrt{ab}J_{m,n}\\ \sqrt{ab}J_{n,m}&bJ_{n,n}},
    \quad
    \pmat{cJ_{m,m}&\sqrt{cd}J_{m,n}\\ \sqrt{cd}J_{n,m}&dJ_{n,n}}
\]
with $a$, $b$, $c$, $d$ to be determined, along with the eigenvalues $\mu_1$ and $\mu_2$.

To determine $\mu_1$ and $\mu_2$, we note that the partition
of $V(X+Y)$ with two cells $V(X)$ and $V(Y)$ is equitable, with quotient
\[
    \pmat{k&n\\ m&\ell}.
\]
Hence $\mu_1$ and $\mu_2$ are the zeros of the characteristic polynomial of this matrix
\[
    x^2 -(k+\ell)x + k\ell-mn,
\]
that is they are equal to
\[
    \frac12 (k+\ell \pm\sqrt{(k-\ell)^2+4mn}).
\]

Now let $u$ and $v$ be distinct vertices in $X$. We determine conditions
for perfect state transfer from $u$ to $v$ in $X+Y$. 
Since $(\hat{F_s})_{u,v}=0$ and $(\hat{E}_r)=(E_r)_{u,v}$,
\begin{align*}
    (H_{X+Y}(t))_{u,v} &= \exp(i\mu_1t)\, (N_1)_{u,v} +\exp(i\mu_2t)\, (N_2)_{u,v} 
        +\sum_{r>1} \exp(i\th_rt)\, (E_r)_{u,v}\\
        &= a\exp(i\mu_1t) +c\exp(i\mu_2t) +\sum_{r>1} \exp(i\th_rt)\, (E_r)_{u,v}.
\end{align*}
On the other hand
\[
    H_X(t)_{u,v} = \frac1m\exp(ikt) +\sum_{r>1} \exp(i\th_rt)\, (E_r)_{u,v}.    
\]
and thus if we have perfect state transfer from $u$ to $v$ in $X$ at time $\tau$,
we will have perfect state transfer between the same vertices in $X+Y$ at time
$\tau$ if
\begin{equation}
    \label{eq:kmumu}
    \frac1m\exp(ik\tau) = a \exp(i\mu_1\tau) + c \exp(i\mu_2\tau).
\end{equation}
As
\[
    I = N_1+N_2 +\sum_{r>1}\hat{E}_r +\sum_{s>1}\hat{F}_s
\]
we see that $a+c=1/m$. Since $\exp(ikt)$, $\exp(i\mu_1\tau)$ and $\exp(i\mu_2\tau)$
are roots of unity, we see that \eqref{eq:kmumu} can hold if and only if
\[
    \exp(ik\tau) = \exp(i\mu_1\tau) =\exp(i\mu_2\tau).
\]
For this we need both $(k-mu_1)\tau$ and $(k-\mu_2)\tau$ to be integer multiples of
$2\pi$, and hence that
\[
    \frac{k-\mu_1}{k-\mu_2} \in \rats.
\]
This can only happen if $(k-\ell)^2+4mn$ is a perfect square.

In our treatment here we have followed Angeles-Canul et al~\cite{Angeles-Canul2009}.
Using these ideas they prove the following results.

\begin{lemma}
    Suppose $Y$ is a $k$-regular graph on $n$ vertices. Then there is perfect state 
    transfer in $\comp{K_2}+Y$ between the vertices of $\comp{K_2}$ if
    \begin{enumerate}[(a)]
        \item
        $\De=\sqrt{k^2+8n}$ is an integer.
        \item
        $n$ is even and $4|k$.
        \item
        The largest power of two that divides $k$ is not equal to the largest power
        that divides $\De$.\qed
    \end{enumerate}
\end{lemma}

\begin{lemma}
    Suppose $Y$ is a $k$-regular graph on $n$ vertices. Then there is perfect state 
    transfer in $K_2+Y$ between the vertices of $K_2$ if
    \begin{enumerate}[(a)]
        \item
        $\De=\sqrt{(k-1)^2+8n}$ is an integer.
        \item
        Both $k-1$ and $n$ are divisible by 8.\qed
    \end{enumerate}
\end{lemma}

The join of $X$ and $Y$ is the complement of the disjoint union of $\comp{X}$ and 
$\comp{Y}$. Hence for regular graphs we can derive results about joins from
information about complements.

\begin{lemma}
    Suppose $X$ is regular graph on $n$ vertices with perfect state transfer from
    $u$ to $v$ at $\tau$. If $\tau$ is an integer multiple of $2\pi/n$, then there
    is perfect state transfer from $u$ to $v$ at time $\tau$ in $\comp{X}$.
\end{lemma}

\proof
Since $X$ is regular, $A$ and $J-I$ commute and so 
\[
    H_{\comp{X}}(t) =\exp(it(J-I)) H_X(-t).
\]
Using the spectral decomposition of $J-I$ we find that
\[
    \exp(it(J-I)) = \exp(i(n-1)t)\frac1n J +\exp(-it)(I-\frac1n J)
\]
and this is a multiple of $I$ if $\exp(int)=1$, that is, if $t$ is an integer
multiple of $2\pi/n$.\qed

Using this lemma, it is immediate that we have perfect state transfer on
$\comp{nK_2}$ (when $n\ge2$) and $\comp{nC_4}$.

\section{The Direct Product}

If $X$ and $Y$ are graphs then their \textsl{direct product} $X\times Y$
is the graph with adjacency matrix 
\[
    A(X) \otimes A(Y).
\]

\begin{lemma}
    \label{lem:dph}
    Suppose $X$ and $Y$ are graphs with respective adjacency matrices $A$
    and $B$ and suppose $A$ has spectral decomposition
    \[
        A = \sum_r \th_r E_r.
    \]
    Then
    \[
        H_{X\times Y}(t) = \sum_r E_r\otimes H_Y(\th_r t).
    \]
\end{lemma}

\proof
First,
\[
    A \otimes B = \sum_r \th_r E_r\otimes B
\]
and since the matrices $E_r\otimes B$ commute,
\[
    H_{X\times Y}(t) = \prod_r \exp(i\th_r E_r\otimes B).
\]
If $E^2=E$ then
\[
    \exp(E\otimes M) = I+ \sum_{k\ge1} \frac1{k!} E\otimes M^k
        = (I-E)\otimes I + E\otimes \exp(M)
\]
and accordingly
\[
    H_{X\times Y}(t) = \prod_r \bigl((I-E_r)\otimes I + E_r\otimes H_Y(\th_r t)\bigr).
\]
Since $E_rE_s=0$ if $r\ne s$ and $\prod_r(I-E_r)=0$, the lemma follows.\qed

\begin{lemma}
    Suppose that $H_Y(\tau)\vc{u}=\ga\vc{v}$ where $\ga=\exp(i\vphi)$
    and $H_Y(2\tau)=\ga^2I$. If the
    eigenvalues $\seq\th1m$ of $X$ are odd integers, then
    \[
        H_{X\times Y}(\tau) = H_X(\vphi)\otimes \ga^{-1}H_Y(\tau).
    \]
\end{lemma}

\proof
Assume that $H_Y(t)=\sum_r\exp(it\th_r)E_r$. If $\th_r$ is an odd integer then
\[
    H_Y(\th_r\tau) = H_Y(2\tau)^{(\th_r-1)/2} H_Y(\tau)
        = \ga^{\th_r} \ga^{-1}H_Y(\tau)
\]
and accordingly
\[
    H_{X\times Y}(t) = \left(\sum_r \ga^{\th_r}E_r\right)\otimes \ga^{-1}H_Y(\tau).
\]
If $\ga=\exp(i\varphi)$ it follows that
\[
    H_{X\times Y}(t) = \left(\sum_r \exp(i\vphi\th_r)E_r\right)\otimes \ga^{-1}H_Y(\tau)
        = H_X(\vphi)\otimes \ga^{-1}H_Y(\tau)
\]
as required.\qed

A closely related result appears as Proposition~2 in Ge et al.~\cite{Ge2010}.
We present two examples provided there.

If the eigenvalues of $X$ are odd integers then $H_X(0)=I$ and $H_X(\pi)=-I$.
The eigenvalues of the $d$-cube are the integers $d-2r$ for $r=0,\ldots,d$, and
it has perfect state transfer at $\pi/2$, with $\ga=i^d$. Hence if $X$ is a graph
with odd integer eigenvalues, then we have perfect state transfer on the
product $X\times Q_{d}$ when $d$ is even. For a second example, 
if $R_d$ denotes the $d$-th Cartesian
power of $P_3$ then $R_d$ has perfect state transfer at time $\pi/\sqrt{2}$
with $\ga=(-1)^d$. Therefore $X\times R_d$ has perfect state transfer at $\pi/\sqrt2$.

Ge et al.~\cite{Ge2010} also give results for the lexicographic product.

\section{Mixing}

Questions about perfect state transfer might be viewed as asking at what times
$t$ does the transition matrix satisfy certain restrictions on its entries.
There are a number of interesting questions of this form.

\subsection{Perfect Mixing}

For the first, we can ask if there is a time $t$ such that all entries
of $H(t)$ have the same absolute value. We say a unitary matrix is \textsl{flat}
if all its entries have the same absolute value and we say that
\textsl{perfect mixing} occurs at time $t$ is $H(t)$ is flat. We have
\[
    H_{K_2}(\pi/4) =\frac1{\sqrt2}\pmat{1&i\\ i&1}
\]
which is flat. Since
\[
    H_{Q_d}(t) =  H_{K_2}(t)^{\otimes d}
\]
we have perfect mixing at time $\pi/4$ on the $d$-cube. Also
\[
    H_{K_4}(t) = \exp(3it)\frac14 J +\exp(-it)\left(I-\frac14J\right)
        = \exp(-it)\left(\exp(4it)\frac14J+I-\frac14J\right)
\]
and thus $H_{K_4}(\pi/4)$ is flat. Consequently any Cartesian product of copies
of $K_2$ and $K_4$ is perfect mixing at time $\pi/4$.

The graphs $K_2$ and $K_4$ are the first two members of a series of graphs: folded cubes.
The \textsl{folded $(d+1)$-cube} is the graph we get from the $d$-cube by joining
each vertex to the unique vertex at distance $d$ from itself. It can also be viewed
as the quotient of the $(d+1)$-cube over the equitable partition formed by the pairs
of vertices at distance $d+1$, which is the origin of the term `folding'.
The first interesting example is the folded 5-cube, often known as the Clebsch graph.
In \cite{Best2008} Best et al.~prove (in our terms) that when $d$ is odd, the folded 
$d$-cube has perfect mixing. 

Ahmadi et al.~\cite{mixcirc2003} prove that $K_3$ is perfect mixing and in \cite{Carlson2006}
Carlson et al.~show that $C_5$ is not. Konno \cite[Section~10.3]{Konno2008} shows 
that $C_6$ is not perfect mixing. We can prove a little more. If $m$ is odd, 
\[
    C_{2m} \cong K_2\times C_m.
\]
Then by \lref{dph}
\[
    H_{K_2\times X}(t) = \frac12\pmat{1&1\\1&1}\otimes H_X(t) 
        + \frac12\pmat{1&-1\\-1&1}\otimes H_X(-t).
\]
If $H_{K_2\times X}(t)$ is flat and
\[
    a := (H_X(t))_{u,v},\quad b :=(H_X(-t))_{u,v}.
\]
then $|a+b|=|a-b|$ and this can hold if and only if $b=ia$. Hence $H(-t)=i H_X(t)$
and
\[
    H_{K_2\times X}(t) = \frac12\pmat{1+i&1-i\\ 1-i&1+i}\otimes H_X(t).
\]
Hence $K_2\times X$ is uniform mixing if and only if $H_X(2t)=-iI$ and $X$ is
uniform mixing at $t$. As $K_3$ is perfect mixing when $t=4\pi/9$, it follows that
$C_6$ is not perfect mixing. Since $C_5$ is not perfect mixing, neither is $C_{10}$.

If perfect mixing occurs at time $\tau$, then $H(\tau)$ is a flat unitary matrix.
Such matrices form an important class of so-called \textsl{type-II} matrices. For further
information see \cite{Chan2007}.

\subsection{Average Uniform Mixing}

For all $t$, each row of the Schur product
\[
    H(t)\circ H(-t)
\]
is a probability density; it is equal to
\[
    \sum_{r,s} \exp(it(\th_r-\th_s))\, E_r\circ E_s
        = \sum_r E_r\circ E_r + 2\sum_{r<s} \cos(t(\th_r-\th_s))\,E_r\circ E_s
\]
and its average value over time is 
\[
    \sum_r E_r\circ E_r.
\]

\begin{lemma}\label{lem:}
    If the $n\times n$ matrices $\seq F1k$ are positive semidefinite and
    $\sum_r F_r$ is a multiple of $J$, then $F_r$ is a multiple of $J$ for each $r$.
\end{lemma}

\proof
If $F$ is positive semidefinite $z=\vc{u}-\vc{v}$, then
\[
    0 \le z^TFz = F_{u,u} + F_{v,v} -2F_{u,v}
\]
and if equality holds
\[
    F_{u,u} = F_{v,v} = F_{u,v}.
\]
Now 
\[
    z^T\left(\sum_r F_r\right)z 
        = \sum_r \bigl((F_r)_{u,u} + (F_r)_{v,v} 
            -2(F_r)_{u,v}\bigr)
\]
Since $F_r$ is positive semidefinite each summand above is non-negative. 
If $\sum_r F_r$ is a multiple of $J$ then the left side is zero, 
and so each summand on the right is zero.\qed

The continuous quantum walk on $X$ is \textsl{average uniform mixing} 
if the average value of $H(t)\circ H(-t)$ is a multiple of $J$. Our next
result is new.

\begin{lemma}
    If $|V(X)|\ge3$ then the continuous quantum walk on $X$ is not average
    uniform mixing.
\end{lemma}

\proof
If $X$ is not connected, it cannot be average uniform mixing, so we assume that $X$
is connected. If the continuous walk on $X$ is average uniform mixing, then
\[
   \sum_r E_r\circ E_r
\]
is a multiple of $J$. Since $E_r$ is positive semidefinite, $E_r\circ E_r$
is positive semidefinite for each $r$ and therefore the previous lemma implies
that $E_r\circ E_r$ is a multiple of $J$.

Since $\th_1$ is the spectral radius of $A$, then entries of $E_1$ are non-negative
and as $E_1\circ E_1$ is a multiple of $J$, it follows that $E_1$ is a multiple of $J$.
Hence we may assume $X$ is $k$-regular. If $n=|V(X)|$ then $E_1=\frac1n J$.
If $r\ne1$ then $E_r$ is flat and since $E_1E_r=0$ we see that 
$n$ is even and exactly half the entries
in any row or column of $E_r$ are negative. This implies that each eigenvalue of $X$
is an integer, with the same parity as the valency $k$. If $|(E_r)_{u,v}| = e$, then
\[
    e = (E_r)_{u,u} = (E_r^2)_{u,u} = ne^2
\]
and thus $e=1/n$. Therefore $\tr(E_r)=1$, which tells us that each eigenvalue of $X$
is simple. Since all diagonal entries of $E_r$ must be positive, they are all equal
and by \cite[Theorem~4.1]{Godsil198051} it follows that $X$ is walk-regular.
From \cite[Theorem~4.8]{Godsil198051} we know that the only connected 
walk-regular graph with simple integer eigenvalues is $K_2$.\qed

Adamczak et al.~\cite{Adamczak2008} prove the above theorem for Cayley graphs of
abelian groups. Our proof follows theirs closely. 

We do not seem to know very much about the average value of $H(t)$. The following
indicates that there may be some surprises.

\begin{theorem}\label{thm:pavg}
    Suppose $\seq E1n$ are the idempotents for the path $P_n$ and let $T$ be the
    permutation matrix such that $T\vc{u}=\vc{n+1-u}$ for $u=1,\ldots,n$. Then
    \[
        \sm r1n E_r\circ E_r = \frac1{2n+2}(2J+I+T).\qed
    \]
\end{theorem}

For a proof and further information, see \cite{Godsil2011}.

\section{Weighted Adjacency Matrices}

Suppose $|V(X)|=n$. We say a symmetric matrix $M$ is a \textsl{weighted adjacency matrix}
for $X$ if $M_{u,v}=0$ for each pair of distinct nonadjacent vertices $u$ and $v$.
So if $M$ is a weighted adjacency matrix for $X$, it is a weighted adjacency matrix
for any graph $Y$ such that $X$ is a subgraph of $Y$. If the off-diagonal entries of
$M$ are $0$ or $\pm1$ and the diagonal entries are zero, we will call $M$ a
\textsl{signed adjacency matrix}. Much of the theory of perfect
state transfer extends to weighted adjacency matrices with very little effort,
since spectral decomposition still applies. 
If $\pi$ is an equitable partition of $X$ then, as we saw, perfect state transfer
on $X$ implies perfect state transfer on the quotient $X/\pi$; the adjacency matrix
of this is weighted and thus we see that information about state transfer on
weighted graphs may have a bearing on the unweighted case. See Ge 
et al.~\cite[Section~5]{Ge2010} for an illustration.
Since the $d$-cube is distance regular, the distance partition relative to any
vertex is equitable and the corresponding quotient is a weighted path. 
As shown by Christandl et al.~\cite{Christandl2005} it follows that for each
$d$ there is a weighting of the edges of a path of length $d$ that admits
perfect state transfer between the end-vertices.
If $\De$ denotes the diagonal matrix of valencies of $X$, then it seems reasonable
to consider the Laplacian $\De-A$ and perhaps the unsigned Laplacian $\De+A$. 
Bose et al.~consider perfect state transfer relative
to the Laplacian in \cite{Bose2008}. (Whether the adjacency matrix or the Laplacian 
is used by a physicist depends on the type of spin interaction postulated.)

To decide which weightings are natural we need to know the physical situation
being modelled, and currently this is largely a matter of speculation. We note
that the fundamental papers \cite{Christandl2005, Saxena2007} focus 
on the unweighted case.

\section{Some Physics}

The states of a quantum system are the 1-dimensional subspaces of a complex
vector space---equivalently the points of a complex projective space.
There are two ways we may avoid admitting that projective geometry is involved.

The first way is to regard vectors non-zero $x$ and $y$ as equivalent if they span the 
same subspace; thus we represent a projective point by an equivalence class of complex
vectors. This is the traditional approach in introductions to quantum physics. 
Here a reversible change of state is modelled by
the application of a unitary operator: if our state is $x$ then the new state
is the subspace spanned by $Ux$, where $U$ is unitary. It is traditional to use unit 
vectors to represent states (which reduces the size of our equivalence classes), and so
then we might say that the new state $y$ is equal to $\ga Ux$. Here $\ga$ is a complex
number of norm 1; physicists call it a phase factor.

The second way is to represent the 1-dimensional space spanned by the nonzero
vector $x$ using the projection
\[
    \frac1{x^*x}xx^*.
\]
If $y= ax$ (and $a\ne0$) then
\[
    \frac1{y^*y}yy^* =\frac1{x^*x}xx^*
\]
If $|x|=1$ and $U$ is unitary then $|Ux|=1$, and the projection on $y$ is
\[
    Uxx^*U^*.
\]
So our phase factors are gone. The map that sends $xx^*$ to $Uxx^*U^*$ is a 
linear map on the space of Hermitian matrices. A quantum state corresponds to
a Hermitian matrix with rank 1 and trace 1. Our unitary operator $U$ is now
\[
    U\otimes U^*.
\]
(Physicists refer to linear maps on spaces of operators as superoperators.)

Note that a positive semidefinite Hermitian matrix can be written as a
sum of matrices of the form $xx^*$ and, if its trace is 1, it can be written
as a convex combination or Hermitian matrices with rank and trace 1. Physicists
refer to the latter as \textsl{pure states} and to a positive semidefinite matrix
with trace 1 as a \textsl{density matrix}.

We turn to our continuous quantum walks, where our operators are the operators $H(t)$.
In the density matrix approach these become
\[
    H(t)\otimes H(t)^* = H(t)\otimes H(-t).
\]
If $H(t)=\exp(itA)$, then
\[
    H(t)\otimes H(-t) = \exp(it(A\otimes I - I\otimes A)).
\]
If $H(t)\vc{u} = \ga\vc{v}$ then $H(-t)\vc{v}= \ga^{-1}\vc{u}$
and so
\[
    H(t)\otimes H(-t)\> \vc{u}\otimes\vc{v} = \vc{v}\otimes\vc{u}.
\]
We also have 
\[
    H(t)\otimes H(-t)\> \vc{u}\otimes\vc{u} = \vc{v}\otimes\vc{v}.
\]
This shows that questions about perfect state transfer on graphs can be translated
to questions about ``phase-free'' perfect state transfer on signed graphs, because
we can view $A\otimes I - I\otimes A$ as the adjacency matrix of a signed version
of the Cartesian square $X\cprod X$.

There is a very interesting recent paper by Pemberton-Ross and Kay 
\cite{Pemberton-Ross2010} using signed adjacency matrices to obtain perfect state
transfer between vertices at distance $n$ in graphs with $cn$ edges 
(for some constant $c$).

\section{Questions}

We list some questions which seem interesting. Unless explicitly stated otherwise,
we consider only unweighted graphs.

\begin{enumerate}[(1)]
    \item
    Is there a graph where we have perfect state transfer from $u$ to $v$,
    but there is no automorphism of $X$ which swaps $u$ and $v$?
    \item
    Let $P_n(a)$ be the path of length $n$ with loops of weight $a$ on each
    end-vertex. Is it true that for each $n$ there is a weight $a$ such that we
    have perfect state transfer between the end-vertices? Casaccino 
    et al.~state in \cite{Casaccino2009} that they have numerical evidence that the 
	answer is yes.
    \item
    Is there some useful theory about the case where $H(t)_{u,u}=0$
    for some time $t$ and some vertex $u$? Or at least when $H(t)\circ I=0$?
    \item
    There are cubelike graphs with perfect state transfer at times $\pi/2$
    and cubelike graphs with perfect state transfer at times $\pi/4$. 
    Are there cubelike graphs with perfect state transfer at time $\tau$, where
    $\tau$ is arbitrarily small?
    \item
    Are there any trees, $K_2$ and $P_3$ aside, on which perfect state transfer
    occurs? [I do not see this as being useful, but it might be fun.]
\end{enumerate}

\section*{Acknowledgements}

I would like to thank the following people, who either provided useful comments on 
early versions of this paper, or helped me understand the material presented in it: 
XiaoXia Fan, Alastair Kay, Dave Morris, Simone Severini, Murray Smith, Christino Tamon, 
Dragan Stevanovi\'c.


\begin{thebibliography}{10}

\bibitem{Adamczak2008}
{\sc W.~Adamczak, K.~Andrew, L.~Bergen, D.~Ethier, P.~Hernberg, J.~Lin, and
  C.~Tamon}, {\em {Non-uniform mixing of quantum walk on cycles}},
  International Journal of Quantum Information, 5 (2007), p.~12. [arxiv:0708.2096]

\bibitem{mixcirc2003}
{\sc A.~Ahmadi, R.~Belk, C.~Tamon, and C.~Wendler}, {\em {On mixing in
  continuous-time quantum walks on some circulant graphs}}, Quantum Information
  and Computation, 3 (2003), pp.~611--618. [arXiv:quant-ph/0209106]

\bibitem{Angeles-Canul2009a}
{\sc R.~J. Angeles-Canul, R.~Norton, M.~Opperman, C.~Paribello, M.~Russell, and
  C.~Tamon}, {\em {On quantum perfect state transfer in weighted join graphs}},
  International Journal of Quantum Information, 7 (2009), p.~16. [arXiv:0909.0431]

\bibitem{Angeles-Canul2009}
\leavevmode\vrule height 2pt depth -1.6pt width 23pt, {\em {Perfect state
  transfer, integral circulants and join of graphs}}, Quantum Information and
  Computation, 10 (2010), pp.~325--342. [arxiv:0907.2148]

\bibitem{Basic2011}
{\sc M.~Ba\v{s}i\'{c}}, {\em {Characterization of circulant graphs having
  perfect state transfer}}, arxiv:1104.1825 (2011), p.~14.

\bibitem{MR2561744}
{\sc M.~Ba\v{s}i\'{c} and M.~D. Petkovi\'{c}}, {\em {Some classes of integral
  circulant graphs either allowing or not allowing perfect state transfer}},
  Appl. Math. Lett., 22 (2009), pp.~1609--1615.

\bibitem{MR2645074}
\leavevmode\vrule height 2pt depth -1.6pt width 23pt, {\em {Perfect state
  transfer in integral circulant graphs of non-square-free order}}, Linear
  Algebra Appl., 433 (2010), pp.~149--163.

\bibitem{MR2523011}
{\sc M.~Ba\v{s}i\'{c}, M.~D. Petkovi\'{c}, and D.~Stevanovi\'{c}}, {\em
  {Perfect state transfer in integral circulant graphs}}, Appl. Math. Lett., 22
  (2009), pp.~1117--1121.

\bibitem{Bernasconi2008}
{\sc A.~Bernasconi, C.~Godsil, and S.~Severini}, {\em {Quantum networks on
  cubelike graphs}}, Physical Review A, 78 (2008), p.~5. [arxiv:0808.0510]

\bibitem{Best2008}
{\sc A.~Best, M.~Kliegl, S.~Mead-Gluchacki, and C.~Tamon}, {\em {Mixing of
  quantum walks on generalized hypercubes}}, International Journal of Quantum
  Information, 6 (2008), pp.~1135--1148. [arxiv:0808.2382]

\bibitem{Bose2003}
{\sc S.~Bose}, {\em {Quantum communication through an unmodulated spin chain}},
  Physical Review Letters, 91 (2003). [arxiv:quant-ph/0212041]

\bibitem{Bose2008}
{\sc S.~Bose, A.~Casaccino, S.~Mancini, and S.~Severini}, {\em {Communication
  in XYZ All-to-All Quantum Networks with a Missing Link}}, International
  Journal of Quantum Information, 7 (2009), pp.~713--723. [arxiv:0808.0748]

\bibitem{Bridges1982}
{\sc W.~G. Bridges and R.~A. Mena}, {\em {Rational G-matrices with rational
  eigenvalues}}, Journal of Combinatorial Theory Series A, 280 (1982),
  pp.~264--280.

\bibitem{Carlson2006}
{\sc W.~Carlson, A.~Ford, E.~Harris, J.~Rosen, C.~Tamon, and K.~Wrobel}, {\em
  {Universal mixing of quantum walk on graphs}}, Quantum Information and
  Computation, 7 (2007), pp.~738--751. [arxiv:quant-ph/0608044]

\bibitem{Casaccino2009}
{\sc A.~Casaccino, S.~Lloyd, S.~Mancini, and S.~Severini}, {\em {Quantum state
  transfer through a qubit network with energy shifts and fluctuations}},
  International Journal of Quantum Information, 7 (2009), pp.~1417--1427.
  [arxiv:0904.4510]

\bibitem{Chan2007}
{\sc A.~Chan and C.~Godsil}, {\em {Type-II matrices and combinatorial
  structures}}, Combinatorica, 30 (2010), pp.~1--24. [arxiv:0707.1836]

\bibitem{Cheung2010}
{\sc W.-C. Cheung and C.~Godsil}, {\em {Perfect state transfer in cubelike
  graphs}}, Linear Algebra and its Applications, 435 (2011), 2468--2474,
  [arxiv:1010.4721]

\bibitem{Christandl2005}
{\sc M.~Christandl, N.~Datta, T.~Dorlas, A.~Ekert, A.~Kay, and A.~Landahl},
  {\em {Perfect transfer of arbitrary states in quantum spin networks}},
  Physical Review A, 71 (2005), p.~12. [arxiv:quant-ph/0411020]

\bibitem{Ge2010}
{\sc Y.~Ge, B.~Greenberg, O.~Perez, and C.~Tamon}, {\em {Perfect state
  transfer, graph products and equitable partitions}}. [arxiv:1009.1340]

\bibitem{Godsil2010b}
{\sc C.~Godsil}, {\em {Controllable subsets in graphs}}, arxiv:1010.3231 (2010), p.~14.

\bibitem{Godsil2010}
\leavevmode\vrule height 2pt depth -1.6pt width 23pt, {\em {When can perfect
  state transfer occur?}}, arxiv:1011.0231 (2010), p.~15,

\bibitem{Godsil2011}
\leavevmode\vrule height 2pt depth -1.6pt width 23pt, {\em {Average mixing of
  continuous quantum walks}}, arXiv:1103.2578 (2011), p.~20, .

\bibitem{Godsil2008}
\leavevmode\vrule height 2pt depth -1.6pt width 23pt, {\em {Periodic Graphs}},
  Electronic J. Combinatorics, 18 (2011), \#23. [arxiv:0806.2074]

\bibitem{cggrbk}
{\sc C.~Godsil and G.~Royle}, {\em {Algebraic Graph Theory}}, Springer, New
  York, 2001.

\bibitem{Godsil2010a}
{\sc C.~Godsil and S.~Severini}, {\em {Control by quantum dynamics on graphs}},
  Physical Review A, 81 (2010), p.~5. [arxiv:0910.5397]

\bibitem{Godsil1987a}
{\sc C.~Godsil and J.~Shawe-Taylor}, {\em {Distance-regularised graphs are
  distance-regular or distance-biregular}}, Journal of Combinatorial Theory,
  Series B, 43 (1987), pp.~14--24.

\bibitem{Godsil198051}
{\sc C.~D. Godsil and B.~D. McKay}, {\em {Feasibility conditions for the
  existence of walk-regular graphs}}, Linear Algebra and its Applications, 30
  (1980), pp.~51--61.

\bibitem{Imrich2000}
{\sc W.~Imrich and S.~Klavzar}, {\em {Product Graphs: Structure and
  Recognition}}, Wiley, 2000.

\bibitem{Kay2010}
{\sc A.~Kay}, {\em {Perfect, efficient, state transfer and its application as a
  constructive tool}}, International Journal of Quantum Information, 08 (2010),
  p.~641. [arxiv:0903.4274]

\bibitem{Kay2011}
\leavevmode\vrule height 2pt depth -1.6pt width 23pt, {\em {The basics of
  perfect communication through quantum networks}}, arxiv:1102.2338 (2011), p.~8.

\bibitem{Kempe2003}
{\sc J.~Kempe}, {\em {Quantum random walks - an introductory overview}},
  Contemporary Physics, 44 (2003), p.~20. [arxiv:0303081]

\bibitem{Kendon2003}
{\sc V.~Kendon}, {\em {Quantum walks on general graphs}}, 
International Journal on Quantum Computation, 4(5) (2006), 791--805.
[arxiv:quant-ph/0306140]

\bibitem{KendTam2011}
{\sc V.~M. Kendon and C.~Tamon}, {\em {Perfect state transfer in quantum walks
  on graphs}}, Journal of Computational and Theoretical Nanoscience, 8 (2011),
  pp.~422--433.

\bibitem{Konno2008}
{\sc N.~Konno}, {\em {Quantum walks}}, in Quantum potential theory, vol.~1954
  of Lecture Notes in Math., Springer, Berlin, 2008, pp.~309--452.

\bibitem{Pemberton-Ross2010}
{\sc P.~J. Pemberton-Ross and A.~Kay}, {\em {Perfect quantum routing in regular
  spin networks}}. Phys. Rev. Lett.\ 106, 020503 (2011), p.~4. [arxiv:1007.2786]

\bibitem{Petkovic2011}
{\sc M.~Petkovi\'c and M.~Ba\v{s}i\'c}, {\em {Further results on the perfect state
  transfer in integral circulant graphs}}, Computers and Mathematics with
  Applications, 61 (2011), pp.~300--312.

\bibitem{Saxena2007}
{\sc N.~Saxena, S.~Severini, and I.~Shparlinski}, {\em {Parameters of integral
  circulant graphs and periodic quantum dynamics}}, 
  International Journal on Quantum Computation, 5(3) (2007), 417--430. [arXiv:quant-ph/0703236]

\bibitem{Stevanovic2011}
{\sc D.~Stevanovi\'c}, {\em {Applications of Graph Spectra in Quantum Physics}},
  in Selected Topics on Applications of Graph Spectra, D.~Cvetkovi\'c and
  I.~Gutman, eds., Belgrade, 2011, Mathematical Institute SANU, pp.~85--111.

\end{thebibliography}

\end{document}